\documentclass{article}
\newtheorem{theorem}{Theorem}
\title{A Compactness Criterion for Real Plane Algebraic Curves}
\author{John Stalker}
\def\mathclap{\mathpalette\mathclapinternal}
\def\mathclapinternal#1#2{%
	\clap{$\mathsurround=0pt#1{#2}$}}
\def\clap#1{\hbox to 0pt{\hss#1\hss}}
\begin{document}
	\maketitle
	\begin{abstract}
		Two sets of conditions are presented for the compactness
		of a real plane algebraic curve, one sufficient and one
		necessary, in terms of the Newton polygon of the defining
		polynomial.
	\end{abstract}
	\section{Introduction}
	We do not, at present, have a complete understanding of the possible
	topologies of real affine plane algebraic curves of given degree.
	Indeed, this is one half of Hilbert's 16th problem.  For any given
	curve the problem is much easier, but is still complicated in comparison
	with the complex case.  A summary of the current state of knowledge,
	and some new results, may be found in~\cite{MR2003h:14088}.
	
	This paper is concerned exclusively with compactness of such curves.
	It is not hard to see that a real affine plane algebraic curve is
	compact if and only if its points on the line at infinity are all
	isolated.  There are finitely many such points and there exists an
	effective procedure to check whether they are
	isolated~\cite{MR90c:14001}.
	That procedure is complicated both to describe and perform, however.
	
	This paper presents a simple necessary criterion and a simple
	sufficient condition, both of which can be stated in terms of the
	Newton polygon of the defining polynomial.  Unfortunately any
	criterion which is both necessary and sufficient will not be
	simple.  For almost all curves with a given Newton polygon, however,
	either the necessary condition fails or the sufficient condition
	holds.  Only for a set of curves of codimension at least one do we
	need to use the more complicated machinery of~\cite{MR90c:14001}.

	In this paper
	\begin{equation}\label{eqn1}
		p(x,y) = \sum c_{k,l} x^k y^l
	\end{equation}
	is a polynomial in two variables with real coefficients and $C$
	is the curve
	\begin{equation}\label{eqn2}
		C = \{ (x,y) \in \mathbf{R}^2 \colon p(x,y) = 0 \},
	\end{equation}
	considered as a subset of $\mathbf{R}^2$ in the usual topology.
	The Newton polygon of $p$ is defined to be the convex hull of the
	set
	\begin{equation}\label{eqn3}
		N = \{ (k,l) \in \mathbf{Z}^2 \colon c_{k,l} \neq 0 \}.
	\end{equation}
	If $E$ is an (oriented) edge of the Newton polygon with endpoints
	$k'_E,l'_E$ and $k''_E,l''_E$ then numbers $d_E$, $p_E$,
	and $q_E$ are defined by
	\begin{equation}\label{eqn4}
		d_E = \gcd(k''_E - k'_E, l''_E - l'_E),
		\qquad
		p_E = \frac{k''_E - k'_E}{d_E},
		\qquad
		q_E = \frac{l''_E - l'_E}{d_E},
	\end{equation}
	and the edge polynomial $e_E \in \mathbf{R}[t]$ by
	\begin{equation}\label{eqn5}
		e_E(t) = \sum_{i=0}^{d_E} c_{k'_E+ip_E,l'_E+iq_E} t^i.
	\end{equation}
	An edge is called outer if it maximizes some linear function $ak+bl$
	on the Newton polygon, where at least one of $a$ or $b$ is positive.

	The purpose of this paper is to establish the following two
	theorems
	\begin{theorem}\label{S}
		For the compactness of $C$ it suffices that
		\begin{enumerate}
			\item $p$ is not divisible by $x$ or $y$, and
			\item the edge polynomials corresponding to outer
			edges have no real zeroes.
		\end{enumerate}
	\end{theorem}
	\begin{theorem}\label{N}
		For the compactness of $C$ it is necessary that
		\begin{enumerate}
			\item\label{N1}
			$p$ is not divisible by $x$ or $y$, and
			\item\label{N2}
			the edge polynomials corresponding to outer
			edges have no real zeroes of odd order.
		\end{enumerate}
	\end{theorem}

	It is possible to give compactness criteria which are both necessary
	and sufficient, but these are much more complicated and require a
	knowledge of $c_{k,l}$ for $(k,l)$ in the interior of the Newton
	polygon, as the example
	\begin{equation}\label{eqn6}
		p(x,y) = x^8 - 4x^6y^2 + 6x^4y^4 - 4x^2y^6 + y^8
		+ cx^2y^2 + 1
	\end{equation}
	shows.  Its zero set is easily seen to be compact if and only if
	$c \ge 0$, but the point $(2,2)$ lies inside the Newton polygon.

	The edge polynomials are polynomials in a single variable, so
	the Sturm test~\cite{dickson} can be used to check for real
	zeroes.

	It is, perhaps, of interest that Theorem~\ref{S} was developed
	for the single example
	\begin{equation}\label{eqn7}
		\begin{array}{r@{\:}c@{\:}l}
			p(x,y) &=&
			72x^{14}
			- (576 + 432y^2)x^{13}
			+ (1947 + 3552y^2 + 1152y^4)x^{12}
			\cr && {}
			- (3504 + 11988y^2 + 10464y^4 + 1440y^6)x^{11}
			\cr && {}
			+ (3452 + 20360y^2 + 38762y^4
				+ 15384y^6 + 720y^8)x^{10}
			\cr && {}
			- (1536 + 16456y^2 + 71800y^4 + 66316y^6
				+ 10536y^8)x^9
			\cr && {}
			+ (2040y^2 + 62966y^4 + 143492y^6
				+ 57803y^8 + 2160y^{10})x^8
			\cr && {}
			- (-4608y^2 + 8608y^4 + 153832y^6
				+ 154672y^8 + 21648y^{10})x^7
			\cr && {}
			+ (-20100y^4 + 48272y^6 + 208760y^8
				+ 83120y^{10} + 2760y^{12})x^6
			\cr && {}
			- (-36120y^6 + 104440y^8 + 151552y^{10}
				+ 20824y^{12})x^5
			\cr && {}
			+ (-33769y^8 + 109100y^{10} + 58958y^{12}
				+ 1908y^{14})x^4
			\cr && {}
			- (-17900y^{10} + 62848y^{12} + 11680y^{14})x^3
			\cr && {}
			+ (-5530y^{12} + 20912y^{14} + 944y^{16})x^2
			\cr && {}
			- (-972y^{14} + 3888y^{16})x + (-81y^{16} + 324y^{18})
		\end{array}
	\end{equation}
	which arises in a problem on wave propagation on a singular
	solution of the Einstein-Maxwell equations in general relativity.
	See \cite{ST} for details.

	\section{Sufficiency}
	Assume that Theorem~\ref{S} is false, i.e. that for some $p$ satisfying
	the conditions there is a sequence of points $(x_n,y_n)$ such that
	\begin{equation}\label{eqn8}
		p(x_n,y_n) = 0
	\end{equation}
	with either $x_n$ or $y_n$ unbounded.  Passing to a subsequence we may
	assume that either
	\begin{equation}\label{eqn9}
		|x_n| \le |y_n| \qquad \mathrm{and} \qquad
		\lim_{n\to\infty} |y_n| = \infty
	\end{equation}
	or
	\begin{equation}\label{eqn10}
		|y_n| \le |x_n| \qquad \mathrm{and} \qquad
		\lim_{n\to\infty} |x_n| = \infty.
	\end{equation}
	Without loss of generality we may assume the former.
	It then follows that
	\begin{equation}\label{eqn11}
		a_n = \frac{\log |x_n|}{\log |y_n|} \le 1
	\end{equation}
	Passing again to a subsequence, we may assume that this quantity
	either tends to a finite limit~$a$ or that it tends to
	$-\infty$.

	Supposing that $\lim_{n\to\infty} a_n = - \infty$, let
	\begin{equation}\label{eqn12}
		L = \max_{(k,l)\in N, k = 0} l
	\end{equation}
	and
	\begin{equation}\label{eqn13}
		A = \min_{(k,l)\in N, k\neq 0} \frac{L - l}{k}.
	\end{equation}
	Then maximum in (\ref{eqn12}) is over a non-empty set since
	$x$ does not divide $p$.  The minimum in (\ref{eqn13}) might be
	over an empty set.  If it is not then, since
	$\lim_{n\to\infty} a_n = - \infty$, we may pass to a subsequence
	where
	\begin{equation}\label{eqn14}
		a_n \le A - 1.
	\end{equation}
	Then, for all $(k,l)\in N$ with $k\neq 0$,
	\begin{equation}\label{eqn15}
		\frac{\log |x_n^k y_n^l|}{\log |y_n|}
		= k a_n + l
		\le k a_n + L - k A
		\le L - k
		\le L - 1
	\end{equation}
	so that
	\begin{equation}\label{eqn16}
		|x_n^k y_n^l| \le |y_n|^{L-1}.
	\end{equation}
	Passing to a subsequence where $|y_n|>1$, the same estimate holds
	for $(k,l)\in N-(0,L)$ such that $k=0$.  It now follows from
	the triangle inequality that
	\begin{equation}\label{eqn17}
		|\sum_{\mathclap{(k,l)\in N-(0,L)}} c_{k,l} x_n^k y_n^l|
		\le |y_n|^{L-1} 
		|\sum_{\mathclap{(k,l)\in N-(0,L)}} c_{k,l}|
	\end{equation}
	Once
	\begin{equation}\label{eqn18}
		|y_n| > |\sum_{\mathclap{(k,l)\in N-(0,L)}} c_{k,l}| / |c_{0,L}|
	\end{equation}
	the triangle inequality shows that
	\begin{equation}\label{eqn19}
		|p(x_n,y_n)| \ge |c_{0,L}||y_n|^L
		- |\sum_{\mathclap{(k,l)\in N-(0,L)}} c_{k,l} x_n^k y_n^l| > 0,
	\end{equation}
	contradicting (\ref{eqn8}).

	Suppose, then, that $\lim_{n\to\infty a_n} = a$.
	Define
	\begin{equation}\label{eqn20}
		m = \max_{(k,l)\in N} (ak + l).
	\end{equation}
	Let $M$ be the set where this maximum is taken.  $M$ is either
	a vertex or an outer edge of the Newton polygon.  In either
	case, define
	\begin{equation}\label{eqn21}
		m''' = \max_{(k,l)\in N-M} (ak + l),
		\qquad m' = \frac{2m + m'''}{3},
		\qquad m'' = \frac{m+m''}{3}.
	\end{equation}
	It follows immediately that $m > m' > m'' > m'''$.

	If $M$ is the single vertex $(k_V,l_V)\in N$ then
	\begin{equation}\label{eqn22}
		\lim_{n\to\infty}
		\frac{\log |c_{k_V,l_V} x_n^{k_V} y_n^{l_V}|}{\log |y_n|}
		= a k_V + l_V = m
	\end{equation}
	so
	\begin{equation}\label{eqn23}
		\log |c_{k_V,l_V} x_n^{k_V} y_n^{l_V}| \ge m'\log |y_n|
	\end{equation}
	and hence
	\begin{equation}\label{eqn24}
		|\sum_{\mathclap{(k,l)\in M}} c_{k,l} x_n^k y_n^l | \ge |y_n|^{m'}
	\end{equation}
	for all sufficiently large $n$.

	If, on the other hand, $M$ is the outer edge $E$ with endpoints
	$(k'_E,l'_E)$ and $(k''_E,l''_E)$ then, since $e_E$ has no real
	zeroes, there is an $\epsilon_E > 0$ such that
	\begin{equation}\label{eqn25}
		|e_E(t)| \ge \epsilon_E
	\end{equation}
	for all~$t$.
	Then
	\begin{equation}\label{eqn26}
		|\sum_{\mathclap{(k,l)\in M}} c_{k,l} x_n^k y_n^l|
		= |x_n^{k'_E} y_n^{l'_E} e_E(x_n^{p_E} y_n^{q_E})|
		\ge \epsilon_E |x_n^{k'_E} y_n^{l'_E}|.
	\end{equation}
	It follows that
	\begin{equation}\label{eqn27}
		\lim_{n\to\infty}
		\frac{\log |\sum_{(k,l)\in M} c_{k,l} x_n^k y_n^l|}{\log |y_n|}
		= a k'_E + l'_E = m
	\end{equation}
	and hence that (\ref{eqn24}) holds for sufficiently large~$n$.

	In either case, for $(k,l) \in N - M$,
	\begin{equation}\label{eqn28}
		\lim_{n\to\infty} \frac{\log |c_{k,l} x_n^k y_n^l|}{\log |y_n|}
		= a k + l \le m'''.
	\end{equation}
	For sufficiently large $n$, we then have
	\begin{equation}\label{eqn29}
		|c_{k,l} x_n^k y_n^l| \le |y_n|^{m''}
	\end{equation}
	and hence
	\begin{equation}\label{eqn30}
		|\sum_{\mathclap{(k,l)\in N - M}} c_{k,l} x_n^k y_n^l|
		\le \#(N-M) |y_n|^{m''}.
	\end{equation}
	For $n$ sufficiently large,
	\begin{equation}\label{eqn31}
		|y_n|^{m'} > \#(N-M) |y_n|^{m''}
	\end{equation}
	and hence, by the triangle inequality,
	\begin{equation}\label{eqn32}
		|p(x_n,y_n)| \ge |\sum_{\mathclap{(k,l)\in M}} c_{k,l} x_n^k y_n^l |
		- |\sum_{\mathclap{(k,l)\in N - M}} c_{k,l} x_n^k y_n^l| > 0,
	\end{equation}
	contradicting (\ref{eqn8}).

	\section{Necessity}
	Condition~\ref{N1} is obviously necessary for compactness.
	Suppose condition~\ref{N2} is violated.  Then there is some
	outer edge $E$ with endpoints $(k'_E,l'_E)$ and $(k''_E,l''_E)$
	for which $e_E$ has a zero of odd order.
	Since $e_E(0) = c_{k'_E,l'_E}$ and $(k'_E,l'_E)$ is an extreme
	point of the Newton polygon we know that $e_E(0) \neq 0$.
	On either side of a zero of odd order
	there are points $t_+$ and $t_-$ where
	\begin{equation}\label{eqn33}
		e_E(t_+) > 0 \qquad \mathrm{and} \qquad e_E(t_-) < 0.
	\end{equation}
	We may take $t_+$ and $t_-$ to have the same sign.
	Since $E$ is an outer edge there are $a$ and $b$, at least one of
	which is positive, such that
	\begin{equation}\label{eqn34}
		a k + b l
	\end{equation}
	is maximized for $(k,l)\in E$.  Let $m$ be this maximum and
	\begin{equation}\label{eqn35}
		m''' = \max_{(k,l) \in N - E} (a k + b l)
	\end{equation}
	Since
	\begin{equation}\label{eqn36}
		\gcd(p_E,q_E) = 1
	\end{equation}
	there are integers $r$ and $s$ such that
	\begin{equation}\label{eqn37}
		r p_E + s q_E = 1.
	\end{equation}
	Setting
	\begin{equation}\label{eqn38}
		x_n(t) = t^r n^a \quad \mathrm{and} \quad y_n(t) = t^s n^b
	\end{equation}
	we see that
	\begin{equation}\label{eqn39}
		x_n(t)^k y_n(t)^l = t^{rk+sl}n^{ak+bl}.
	\end{equation}
	If $k = k'_E + i p_E$ and $l = l'_E + i q_E$ then
	\begin{equation}\label{eqn40}
		x_n(t)^k y_n(t)^l = t^{rk'_E + sl'_E + i}n^m.
	\end{equation}
	From this it follows that
	\begin{equation}\label{eqn41}
		\sum_{\mathclap{(k,l) \in E}} c_{k,l} x_n(t)^k y_n(t)^l
		= t^{rk'_E + sl'_E} n^m e_E(t).
	\end{equation}
	If $(k,l) \in N - E$ then
	\begin{equation}\label{eqn42}
		\sum_{\mathclap{(k,l) \in N - E}} c_{k,l} x_n(t_+)^k y_n(t_+)^l
		= \sum_{\mathclap{(k,l) \in N - E}} c_{k,l} t_+^{rk+sl}n^{ak+bl}
	\end{equation}
	and hence, by the triangle inequality,
	\begin{equation}\label{eqn43}
		|\sum_{\mathclap{(k,l) \in N - E}} c_{k,l} x_n(t_+)^k y_n(t_+)^l|
		\le \sum_{\mathclap{(k,l) \in N - E}} c_{k,l} t_+^{rk+sl}| n^{m'''}
	\end{equation}
	If $n$ is sufficiently large that
	\begin{equation}\label{eqn44}
		t^{rk'_E + sl'_E} e_E(t_+) n^m
		> \sum_{\mathclap{(k,l) \in N - E}} c_{k,l} t_+^{rk+sl}| n^{m'''}
	\end{equation}
	then another application of the triangle inequality shows that
	\begin{equation}\label{eqn45}
		p(x_n(t_+),y_n(t_+))
	\end{equation}
	has the same sign as
	\begin{equation}\label{eqn46}
		t_+^{rk'_E + sl'_E}
	\end{equation}
	Similarly, for large $n$,
	\begin{equation}\label{eqn47}
		p(x_n(t_-),y_n(t_-))
	\end{equation}
	has the opposite sign from
	\begin{equation}\label{eqn48}
		t_-^{rk'_E + sl'_E}
	\end{equation}
	and hence $p(x_n(t_+),y_n(t_+))$ and $p(x_n(t_-),y_n(t_-))$ are
	of opposite sign.
	By the intermediate value theorem there is then a $t_n$ between
	$t_-$ and $t_+$ for which
	\begin{equation}\label{eqn49}
		p(x_n(t_n), y_n(t_n)) = 0.
	\end{equation}
	Since at least one of $a$ or $b$ is positive at least one of
	$x_n(t_n)$ or $y_n(t_n)$ is unbounded as $n$ tends to infinity.
	The curve $C$ is therefore not compact.
	\bibliographystyle{plain}
	\bibliography{rag}

\begin{thebibliography}{1}

\bibitem{MR2003h:14088}
Maria~Jesus de~la Puente.
\newblock Real plane algebraic curves.
\newblock {\em Expo. Math.}, 20(4):291--314, 2002.

\bibitem{dickson}
L.E. Dickson.
\newblock {\em First Course in the Theory of Equations}.
\newblock John Wiley \& Sons, New York, 1922.

\bibitem{MR90c:14001}
Dominique Duval.
\newblock Rational {P}uiseux expansions.
\newblock {\em Compositio Math.}, 70(2):119--154, 1989.

\bibitem{ST}
John Stalker and Abdolreza~Shadi Tahvildar-Zadeh.
\newblock Scalar waves on a naked-singularity background.
\newblock 2003.

\end{thebibliography}
\end{document}